\documentclass{amsart}
\usepackage{amssymb}
\newtheorem{thm}{Theorem}
\newtheorem{lem}{Lemma}

\usepackage{comment}

\newcommand{\p}{^{\prime}}
\usepackage{hyperref}

\title{Numbers which are only orders of abelian or nilpotent groups}
\author{Matthew Just}

\begin{document}

\maketitle

\begin{abstract}
    Refining a result of Erd\H{os} and Mays, we give asymptotic series expansions for the functions $A(x)-C(x)$, the count of $n\leq x$ for which every group of order  $n$ is abelian (but not all cyclic), and $N(x)-A(x)$, the count of $n\leq x$ for which every group of order $n$ is nilpotent (but not all abelian).
\end{abstract}

\section{Introduction}

A number $n$ is called \textit{cyclic} if every group of order $n$ is cyclic. Call a number $n$ \textit{abelian} if every group of order $n$ is abelian, and a number $n$ \textit{strictly abelian} if $n$ is abelian and there is at least one group of order $n$ that is not cyclic. Similarly, call $n$ {\it nilpotent} if every group of order $n$ is nilpotent, and {\it strictly nilpotent} if there is at least one group of order $n$ that is not abelian.

Characterizations of cyclic, abelian, and nilpotent numbers can be found in the work of Dickson \cite{Dic05}, Szele \cite{Sze47}, and Pazderski \cite{Paz59}. A number $n$ is cyclic if and only if $(n,\varphi(n))=1$, where $(a,b)$ is the greatest common divisor of $a$ and $b$ and $\varphi(n)$ is Euler's totient function. Define the multiplicative function $\psi(n)$ whose value on a prime power is \[\psi(p^a)=(p^a-1)(p^{a-1}-1)\ldots(p-1).\] Note that $\psi(n)=\varphi(n)$ if and only if $n$ is squarefree. A number $n$ is abelian if and only if $n$ is cubefree and $(n,\psi(n))=1$. Since $(n,\varphi(n))=1$ implies that $n$ is squarefree, a number is strictly abelian if and only if $n$ is cubefree, not squarefree, and $(n,\psi(n))=1$. Similarly, a number $n$ is nilpotent if and only if $(n,\psi(n))=1$, and strictly nilpotent if $n$ is not cubefree.

Let $C(x)$ be the number of cyclic numbers not exceeding $x$, $A(x)$ be the number of abelian numbers not exceeding $x$, and $N(x)$ be the number of nilpotent numbers not exceeding $x$. The number of strictly abelian numbers and strictly nilpotent numbers is then $A(x)-C(x)$ and $N(x)-A(x)$, respectively. Erd\H{o}s \cite{Erd48} showed that 
\[C(x)\sim \frac{e^{-\gamma} x}{\log_3 x}\]
as $x\rightarrow\infty$, where $\log_k x$ is the $k$th iterate of the natural logarithm. Mays \cite{May78} showed that $A(x)\sim C(x)$. This turns out to be a special case of a theorem of Scourfield \cite{Sco76}, who showed that if $S(x)$ counts the number of $n$ with $(n,f(n))=1$ where $f(n)$ is a specific type of multiplicative function (of which $\varphi(n)$ is an example), then $S(x)\sim C(x)$. As $\psi(n)$ also satisfies the conditions outlined by Scourfield for $f(n)$, it follows that $N(x)\sim C(x)$.

Erd\H{o}s and Mays \cite{EM88} subsequently showed that \[A(x)-C(x) \sim \frac{c x}{\log_2 x (\log_3 x)^2}\]
and
\[N(x)-A(x) \sim \frac{c x}{(\log_2 x)^2 (\log_3 x)^2}\]
as $x\rightarrow \infty$, for some constant $c$. This constant is explicitly given as $c=e^{-\gamma}$ in a paper of Narlikar and Srinivasan \cite{NS88}.

Improving on a result of Begunts \cite{Beg01}, Pollack \cite{Pol20} has shown that $C(x)$ admits an asymptotic series expansion (in the sense of Poincar\'e). Specifically,
\[C(x) = \frac{e^{-\gamma}x}{\log_3 x} \left( 1-\frac{\gamma}{\log_3 x} + \frac{\gamma^2+\frac{1}{12}\pi^2}{(\log_3 x)^2} -\frac{\gamma^3+\frac{1}{4}\gamma\pi^2+\frac{2}{3}\zeta(3)}{(\log_3 x)^3}+ \dots \right).\]

Our main results give analogous expansions for $A(x)-C(x)$ and $N(x)-A(x)$. The proofs of our main theorems use ideas similar to those use by Pollack \cite{Pol20}, however the computations are more intricate. Because strictly abelian and strictly nilpotent numbers are not squarefree, we must consider more congruence conditions on the prime divisors. Interestingly, the coefficients of the Gamma function expanded about $w=2$ play a central role in the expansion of $N(x)-A(x)$, while the coefficients of the Gamma function expanded about $w=1$ are relevant for the expansions of $A(x)-C(x)$ and $C(x)$.

\begin{thm}
    For any positive integer $N$, there is a sequence of real numbers $b_0=1$, $b_1$, $b_2$, \ldots, $b_N$ such that
    \[A(x)-C(x) = \frac{xe^{-\gamma}}{\log_2 x (\log_3 x)^2}\left(\sum_{k=0}^N \frac{b_k}{(\log_3 x)^k} \right) + O_N\left( \frac{x}{\log_2 x (\log_3 x)^{N+3}}\right).\]
\end{thm}

The constants $b_k$ are determined as follows. Let $C_k$ be the coefficients of the series expansion of $\Gamma(1+w)$ about $w=0$. Let $c_0=1$, $c_1=-\gamma$, $c_2=\gamma^2+\frac{1}{12}\pi^2$, and in general $c_k$ be determined by the formal relation \[c_0+c_1z+c_2z^2+\ldots = \exp\left(0!C_1z+1!C_2z^2+2!C_3z^3+\ldots \right).\] Then 
\[b_k = \sum_{i+j=k} j!c_{i}C_j.\]
For instance,
\[b_0=1, \ \ \ b_1=-2\gamma, \ \ \ b_2=3\gamma^2+\frac{1}{4}\pi^2, \ \ \ b_3=-4\gamma^3-\gamma\pi^2-\frac{8}{3} \zeta(3).\]

\begin{thm}
    For any positive integer $N$, there is a sequence of real numbers $d_0=1$, $d_1$, $d_2$, \ldots, $d_N$ such that
    \[A(x)-C(x) = \frac{xe^{-\gamma}}{(\log_2 x)^2 (\log_3 x)^2}\left(\sum_{k=0}^N \frac{d_k}{(\log_3 x)^k} \right) + O_N\left( \frac{x}{\log_2 x (\log_3 x)^{N+3}}\right).\]
\end{thm}

The constants $d_k$ are determined as follows. Let $c_k$ be defined as above, and $D_k$ be the coefficients of the series expansion of $\Gamma(2+w)$ about $w=0$. Then
\[d_k=\sum_{i+j=k} j!c_{i}D_j.\]
For instance,

\[d_0=1, \ \ \ d_1=1-2\gamma, \ \ \ d_2=\frac{1}{6}\left(-12\gamma + 15\gamma^2 + \pi^2\right),\]
\[d_3=\frac{1}{6}\left(15\gamma^2 -16\gamma^3+\pi^2-3\gamma \pi^2-6\zeta(3)\right).\]

\section{Lemmata}

The first Lemma is due to Norton \cite{Nor76} and Pomerance \cite{Pom77}. Unless otherwise stated, all implied constants are absolute.

\begin{lem} \label{mert}
    If $a$ and $m$ are two integers such that $m>0$ and $(a,m)=1$, then for $x\geq 3$ \[\sum_{\substack{p\leq x \\ p\equiv a(m)}} \frac{1}{p} = \frac{\log_2 x}{\phi(m)} + \frac{1}{p_{a,m}} +  O\left( \frac{\log 2m}{\phi(m)}\right),\]
    where $p_{a,m}$ is the least prime congruent to $a$ modulo $m$.
\end{lem}

The next result is a special case of the Fundamental lemma of the Selberg sieve, see \cite{HR74}.

\begin{lem}\label{flss}
    For $x\geq y \geq 3$, let $P$ be the product of a set of distinct primes not exceeding $y$. Then the number $n\leq x$ not divisible by any prime dividing $P$ is equal to 
    \[x\prod_{p\mid P} \left(1-\frac{1}{p}\right)\left(1+ O\left(\exp\left(-\frac{1}{2}\frac{\log x}{\log y} \right) \right)\right).\]
\end{lem}

The next result is an application of Lemma \ref{mert} and Lemma \ref{flss}.

\begin{lem}\label{sumcoprime}
    If $a$ and $m$ are two integers such that $m>0$ and $(a,m)=1$, then for $x\geq 3$
    \[\sum_{\substack{n\leq x\\p\mid n \Rightarrow p\not\equiv a(m)}}1 = O\left(x\exp\left(-\frac{\log\log x}{\phi(m)}\right)\right),\]
    where the sum is over the $n\leq x$ that are not divisible by a prime congruent to $a$ modulo $m$.
\end{lem}


The following Lemma is a consequence of the Brun-Titschmarsh inequality \cite{Tit30} and partial summation.

\begin{lem}\label{mertalt}
    If $a$ and $m$ are two integers such that $m>0$ and $(a,m)=1$, then for $x\geq 3$ \[\sum_{\substack{m\leq p\leq x \\ p\equiv a(m)}} \frac{1}{p} = O\left(\frac{\log_2 x}{\phi(m)}\right).\]
\end{lem}

The following two results estimate the sum of the reciprocals of the squares/cubes/fourth powers of the primes greater than some fixed point. They can be deduced by the prime number theorem and partial summation.

\begin{lem}\label{primerecsq}
        The sum of the square of the reciprocal of the primes larger than $T$ as $T\rightarrow \infty$:
        \[\sum_{p>T}\frac{1}{p^2} = \frac{1}{T\log T} (1+o(1)). \]
    \end{lem}
 
 \begin{lem}\label{primereccu}
        The sum of the cube of the reciprocal of the primes larger than $T$ as $T\rightarrow \infty$:
        \[\sum_{p>T}\frac{1}{p^3} = \frac{1}{2T^2\log T} (1+o(1)). \]
    \end{lem}   
\begin{lem}\label{primerec4}
    The sum of the fourth power of the reciprocal of the primes larger than $T$ as $T\rightarrow \infty$:
    \[\sum_{p>T}\frac{1}{p^4} = \frac{1}{3T^3\log T}(1+O(1)).\]
\end{lem}

The final Lemma is a version of Mertens' third theorem incorporating the de la Vall\'ee Poussin error term in the prime number theorem. We will use this error term in many places, and we will use $K$ to represent a (possibly changing) fixed constant.

\begin{lem}\label{mertgood} For $z>1$
    \[\prod_{p\leq z} \left(1-\frac{1}{p} \right) = \frac{e^{-\gamma}}{\log z} (1 + \exp(-K\sqrt{\log z})).\]
\end{lem}

\section{Strictly abelian numbers: proof of Theorem 1}
    
    Let  $y=\log_2 x/\log_3 x$, $z=\log_2 x \exp(\sqrt{\log_3 x})$, and $p$ be a prime in the interval $(y,z]$.  Define $S(x;p,0)$ to be size of the set of $n\leq x$ such that $n=mp^2$, where $m$ is a a product of distinct primes larger than $y$, none of which are equal to $p$.
    
    Furthermore, define $S(x;p,k)$ to be size of the subset of $S(x;p,0)$ with the further restriction that for $n=mp^2$ in $S(x;p,k)$, there are exactly $k$ primes in the interval $(y,z]$ dividing $m$, and at least one of these primes $q$ has the property that there is a prime $z<r \leq x^{1/\log_2 x}$ dividing $m$ that is congruent to 1 modulo $q$ or 1 modulo $p$.
    
    Now define
    \[B(x) = \sum_{y<p\leq z}\left(S(x;p,0)-\sum_{1\leq k \leq \log_3 x} S(x;p,k) \right). \]
    We will show that $A(x)-C(x) = B(x) + O(x/z)$.
    
    Suppose that $n$ is counted by $A(x)-C(x)$, but not $B(x)$. Then one of the following holds:
    
    \begin{enumerate}
        \item There is a prime $q>z$ such that $q^2\mid n$.
        \item There are two distinct primes in the interval $(y,z]$, $q$ and $r$, such that $q^2r^2$ divides $n$.
        \item There is a prime $q$ in the interval $[2,y]$ dividing $n$.
    \end{enumerate}
    
    The number of $n\leq x$ for which (1) holds is at most $\sum_{q>z}x/q^2=O(x/z)$.
    
    For (2), the number of $n\leq x$ with this property is at most
    \begin{align*}
        \sum_{y<q<r\leq z} \frac{x}{q^2 r^2}&\leq x\left(\sum_{q>y} \frac{1}{q^2} \right)^2\\
        &= \frac{x}{y^2(\log y)^2}(1+o(1)),
    \end{align*}
    where we have used Lemma \ref{primerecsq}. This quantity is $O(x/z)$.
    
    Suppose (3) holds. Then either $q^2 \mid n$, or $q^2\nmid n$. Suppose that $q^2\mid n$. Then $n=q^2 m$ where $m$ is not divisible by any prime congruent to 1 modulo $q$. By Lemma \ref{sumcoprime}, the total number of such $n$ is
    \[\ll \sum_{q\leq y} \frac{x}{q^2} \exp\left(-\frac{\log_2 x}{q} \right)\ll x\int_2^y \frac{\exp\left(-\frac{\log_2 x}{t} \right)}{t^2}\, dt = O\left(\frac{x}{(\log_2 x )^2} \right),\] 
    which is $O(x/z)$. Now if $q^2\nmid n$, then $n=p^2 m$ for some $p>y$. Furthermore, for a fixed $p$, we have $m=q \ell\leq x/p^2$ where $\ell$ cannot be divisible by a prime 1 modulo $q$. Again using Lemma \ref{sumcoprime}, the number of $n$ of this form is 
    \[\ll \sum_{p>y} \frac{x}{p^2} \sum_{q\leq y} \frac{1}{q} \exp\left(-\frac{\log_2 x}{q} \right) \ll \sum_{p>y} \frac{x\log_2 y}{p^2 \log_2 x} = O\left( \frac{x\log_2 y}{y \log y \log_2 x} \right),  \]
    which is $O(x/z)$.

    Suppose that $n$ is counted by $B(x)$, but not $A(x)-C(x)$. Then one of the following holds:
    
    \begin{enumerate}
        \item There is a prime $q>z$ such that $q$ divides $n$ and $\phi(n)$.
        \item $m$ has more than $\log_3 x$ prime factors in $(y,z]$.
        \item There is a prime $q$ in $(y,z]$ dividing $n$ such that there is a prime $r>x^{1/\log_2 x}$ dividing $m$ with $r$ congruent to 1 modulo $q$ or 1 modulo $p$.
        \item There is a prime $q$ in the interval $(y,z]$ dividing $m$ such that $p$ is congruent to $\pm 1$ modulo $q$.
        \item There are two distinct primes $r$ and $s$ in $(y,z]$ dividing $n$ such that $s$ is congruent to 1 modulo $r$.
    \end{enumerate}
    
    Suppose that (1) holds. Then there is a prime $r$ that is congruent to 1 modulo $q$ dividing $n$. For a fixed prime $p\in (y,z]$, the number of such $n=p^2 m$ of this form is at most
    \[\sum_{q>z}\sum_{\substack{r\leq x \\ r\equiv 1(q)}}\frac{x}{p^2q r}  \ll \sum_{q>z} \frac{x\log_2 x}{p^2q^2} \ll \frac{x \log_2 x}{p^2 z\log z}\]
    using Lemmas \ref{mertalt} and \ref{primerecsq}, and summing over primes $p$ in the interval $(y,z]$, we see that this quantity is $O(x/z)$.
    
    The number of $n$ for which $(2)$ holds is at most
    \[\sum_{y<p\leq z} \frac{x}{p^2}\sum_{k>\log_3 x}\left(\sum_{y<q\leq z} \frac{1}{q} \right)^k.\]
    The sum $\sum_{y<q\leq z} 1/q= O(\frac{1}{\sqrt{\log_3 x}})$ by Mertens' theorem which, for $x$ sufficiently large, will be less than $1/2$ for all $k>\log_3 x$. Therefore the number of $n$ for which $(2)$ holds is
    \[\ll\sum_{y<p\leq z} \frac{x}{p^2} 2^{-\log_3 x} = O(x/z). \]
    
    The number of $n$ for which $(3)$ holds is \[\ll \sum_{y<p\leq z} \frac{x}{p^2}\sum_{y<q\leq z} \frac{1}{q} \left(\sum_{\substack{x^{1/\log_2 x} < r \leq x \\ r \equiv 1 (q)}} \frac{1}{r}+\sum_{\substack{x^{1/\log_2 x} < r \leq x \\ r \equiv 1 (p)}} \frac{1}{r}\right). \]
    By Lemma \ref{mert}, this is
    \[\ll \sum_{y<p\leq z}\frac{x}{p^2} \sum_{y<q\leq z}\left( \frac{\log_3 x}{q^2} +\frac{\log_3 x}{qp}\right)\ll \frac{x\log_3 x}{ y^2 (\log y)^2}+\frac{x}{y^2 \log y}=O(x/z).\] Here we use that $\sum_{y<q\leq z}\frac{1}{q} \ll \frac{1}{\log_3 x}$ by Mertens' theorem.
    
    For $(4)$, first we count the number of $n$ with a prime divisor $q$ in $(y,z]$ such that $p$ is 1 modulo $q$. The number of such $n$ is at most
    \[\sum_{y<q\leq z} \frac{x}{q} \sum_{\substack{y<p\leq z \\ p\equiv 1 (q) }} \frac{1}{p^2}.  \]
    Now if $p=kq+1$, then $p^2\geq k^2 q^2$. Thus the above quantity is
    \[\ll \sum_{y<q\leq z} \frac{x}{q}  \sum_{m\geq 1} \frac{1}{m^2q^2} \ll \sum_{y<q\leq z} \frac{x}{q^3} =O(x/z),\] by Lemma \ref{primereccu}. The case when $p$ is $-1$ modulo $q$ is similar.
    
    For $(5)$, let $r$ be a prime in $(y,z]$ dividing $n$. Then any primes $s$ in $(y,z]$ dividing $n$ congruent to 1 modulo $r$ will be of the form $s=kr+1$. We consider two cases, $r=p$ and $r\neq p$. The case when $s=p$ is covered in $(4)$, so we may assume $s\neq p$. In the case that $r=p$, the number of $n\leq x$ such that $(5)$ holds is at most
    \[\sum_{y<p<z}\sum_{\substack{y<s<z\\s=kp+1}}\frac{x}{p^2s}\ll\sum_{y<p\leq z} \sum_{k\leq z/y} \frac{x}{kp^3} = O\left( \frac{x\log(z/y)}{y^2\log y} \right) = O(x/z).  \]
    In the case that $r\neq p$, the number of $n\leq x$ such that $(5)$ holds is at most
    \[\sum_{y<p<z}\sum_{y<r<z}\sum_{\substack{y<s<z \\ s=kr+1}}\frac{x}{p^2sr}\ll\sum_{y<p\leq z} \sum_{y<r\leq p} \sum_{k\leq z/y} \frac{x}{kp^2r^2 } = \frac{x\log(z/y)}{y^2(\log y)^2} = O(x/z).\]
    
    Now we estimate $S(x;p,0)$. For each $n$ counted by $S(x;p,0)$, we have $n=p^2 m$ where $m$ is not divisible by $p$ or any prime smaller than $y$. By Lemma \ref{flss}, the number of such $n$ equals
    \[\frac{x}{p^2}\left(1-\frac{1}{p}\right) \prod_{q\leq y} \left(1-\frac{1}{q}\right)\left(1+O\left(\exp\left( -\frac{1}{2}\frac{\log x}{\log y} \right)\right)\right). \]
    By applying Lemma \ref{mertgood} to the product, we see that this quantity equals
    \[\frac{xe^{-\gamma}}{p^2\log y} + O\left(\frac{x}{p^2\exp(K\sqrt{\log_3 x})}\right).\]

    Now to estimate $S(x;p,k)$, where $1\leq k \leq \log_3 x$. For a fixed $k$, write $n=q_0^2 q_1q_2\ldots q_k m$ for a number counted by $S(x;p,k)$, where $q_0=p$, and  $q_1$, $q_2$, \ldots, $q_k$ are increasing distinct primes in $(y,z]$. Furthemore, $m$ is only divisible by primes larger than $z$, and there is a prime $r\leq x^{1/\log_2 x}$ dividing $m$ such that $r$ is congruent to 1 modulo $q_i$ for some $i$.
    
    The number of such $n$ equals
    \begin{align}
    \sum_{\substack{y<q_1<\ldots<q_k\leq z \\ q_i \neq q_0}}\frac{x}{q_0^2q_1\ldots q_k}\prod_{s\leq z}\left( 1-\frac{1}{s} \right)\left(1-\prod_{\substack{z< r \leq x^{1/\log_2 x} \\ r \equiv 1  (q_i) \ \text{for some i}}}\left(1-\frac{1}{r}\right) \right), \label{sk}
    \end{align}
    with an error term (by Lemma \ref{flss}) that is $O(x/p^2\exp(K\sqrt{\log_3 x}))$ .
    The product over primes $r$ equals
    
    \[\prod_{\substack{z< r \leq x^{1/\log_2 x} \\ r \equiv1  (q_i) \ \text{for some i}}}\left(1-\frac{1}{r}\right) = \exp\left( -\sum_{\substack{z< r \leq x^{1/\log_2 x} \\ r \equiv1 \ (q_i)  \text{for some i}} }\frac{1}{r} \right)(1+O(1/z)).\]
    
    The inner sum here equals
    
    \begin{align*}
    \sum_{\substack{z< r \leq x^{1/\log_2 x} \\ r \equiv1 \ (q_i)  \text{for some i}} }\frac{1}{r} &= \sum_{i=0}^k \sum_{\substack{z< r \leq x^{1/\log_2 x} \\ r \equiv1  (q_i) \ \text{for some i}}}\frac{1}{r}+ O\left(\sum_{0\leq i < j \leq k}\sum_{\substack{z< r \leq x^{1/\log_2 x} \\ r \equiv1 (q_i q_j)}} \frac{1}{r} \right).
    \end{align*}
    
    The error term here (using Lemma \ref{mertalt}) is
    
    \[\ll {k+1 \choose 2} \frac{(\log_3 x)^2}{\log_2 x} \ll \frac{(\log_3 x)^4}{\log_2 x}. \]

    The main term equals (Lemma \ref{mert})
    
    \begin{align*}
    \sum_{i=0}^k \sum_{\substack{z< r \leq x^{1/\log_2 x} \\ r \equiv1  (q_i) \ \text{for some i}}}\frac{1}{r} &= \sum_{i=0}^k \frac{\log_2 x}{q_i} + O\left( \frac{(\log_3 x)^3}{\log_2 x} \right).
    \end{align*}
    
    Therefore,
    
   \[ \prod_{\substack{z< r \leq x^{1/\log_2 x} \\ r \equiv1  (q_i)  \ \text{for some i}}}\left(1-\frac{1}{r}\right) = \prod_{i=0}^k \exp\left( \frac{\log_2 x}{q_i} \right) + O\left( \frac{(\log_3 x)^4}{\log_2 x}\right).\]
   
   We can estimate the product of primes up to $z$ by
   \[\prod_{p\leq z} \left(1-\frac{1}{p} \right) = \frac{e^{-\gamma}}{\log z}(1+O(\exp(-K\sqrt{\log _3 x}))\]
   using Lemma \ref{mertgood}, so all together we have that (\ref{sk}) equals
   \[ \sum_{\substack{y<q_1<\ldots<q_k\leq z \\ q_i\neq p}} \left(\frac{e^{-\gamma} x}{\log z} \left( \frac{1}{p^2q_1\ldots q_k} - \frac{1}{p} \prod_{i=0}^k \frac{\exp(-\log_2 x/q_i)}{q_i} \right) + O\left( \frac{x}{p^2q_1\ldots q_k \exp(K\sqrt{\log_3 x})}\right)\right).\]
   The error term is $O(x/(p^2\exp(K\sqrt{\log_3 x})))$.
   
   Let $\sigma = \sum_{y<q\leq z}\frac{1}{q}$. Note that $\sigma \sim \frac{1}{\sqrt{\log_3 x}}$ by Mertens' theorem. We then have the upper bound
   \[\sum_{\substack{y<q_1<\ldots<q_k\leq z \\ q_i\neq p}} \frac{1}{q_1\ldots q_k} \leq \frac{1}{k!} \sigma^k.\] For a lower bound, notice that
   \begin{align*}
   \sum_{\substack{y<q_1,\ldots,q_k\leq z \\ q_i\neq p,\text{ $q_i$ distinct}}} \frac{1}{q_1\ldots q_k} &= \sum_{\substack{y<q_1,\ldots,q_{k-1}\leq z \\ q_i\neq p,\text{ $q_i$ distinct}}} \frac{1}{q_1\ldots q_{k-1}} \sum_{\substack{y<q_k\leq z \\ q_k \neq p\text{ or } q_i \\ \text{ for any $i$}}} \frac{1}{q_k} \\
   &\geq \sum_{\substack{y<q_1,\ldots,q_{k-1}\leq z \\ q_i\neq p,\text{ $q_i$ distinct}}} \frac{1}{q_1\ldots q_{k-1}} \left(\sigma - \frac{k}{y} \right).
   \end{align*}
   
   We can apply a similar lower bound for the remaining $q_i$s to get
   \begin{align*}
   \sum_{\substack{y<q_1<\ldots<q_k\leq z \\ q_i\neq p}} \frac{1}{q_1\ldots q_k} \geq \frac{1}{k!}\prod_{i=1}^k\left(\sigma-\frac{i}{y}\right) &= \frac{1}{k!}\left(\sigma - \frac{(\log_3 x)^2}{\log_2 x} \right)^k \\
    &=\frac{1}{k!}\sigma^k + O\left( \frac{(\log_3 x)^4}{k!\log_2 x} \right),
    \end{align*}
    where we have applied Bernoulli's inequality and used that $\sigma \gg 1/\sqrt{\log_3 x} $ for large $x$ by Mertens' theorem. Therefore,
    \[\sum_{\substack{y<q_1<\ldots<q_k\leq z \\ q_i\neq p}} \frac{1}{q_1\ldots q_k} = \frac{1}{k!}\sigma^k + O\left( \frac{(\log_3 x)^4}{k!\log_2 x} \right).\]
    
    Now set $\tau = \sum_{y<q\leq z} \exp(-\log_2 x/q)/q$. By an analogous argument as above, we have 
    \[\sum_{\substack{y<q_1<\ldots<q_k\leq z \\ q_i\neq p}}\prod_{i=1}^k \frac{\exp(-\log_2 x/q_i)}{q_i} =\frac{1}{k!}\tau^k + O\left( \frac{(\log_3 x)^4}{k!\log_2 x} \right) .\]
    
    Collecting estimates, we have that (\ref{sk}) equals
    \[\frac{e^{-\gamma}x}{k!\log z} \left(\frac{ \sigma^k-\tau^k\exp(-\log_2 x/p)}{p^2}\right)  + O\left( \frac{x(\log_3 x)^3}{p^2 k!\log_2 x} + \frac{x}{p^2\exp( K \sqrt{\log_3 x})} \right).\]
    
    Now
    
    \[\sum_{k\leq \log_3 x} S(x;p,k) = \frac{e^{-\gamma}x}{p^2\log z} \left(\exp(\sigma) - \exp(\tau -(\log_2 x)/p)  \right) + O\left( \frac{x}{p^2\exp(K\sqrt{\log_3 x})} \right).\]
    Therefore, since $\exp(\sigma) = \frac{\log z}{\log y} (1+O(\exp(-K\sqrt{\log_3 x})))$ for large $x$ by Mertens' theorem,
    \[B(x) = \sum_{y<p\leq z}\left(\frac{e^{-\gamma }x}{p^2\log z} \exp\left( \tau-\frac{\log_2 x}{p} \right) + O\left( \frac{x}{p^2\exp(K\sqrt{\log_3 x})}\right)\right).\]
    The error term is $O(x/z)$ after summing on $p$.
    
    A calculation identical to one done by Pollack \cite{Pol20} shows that
    \[\exp(\tau) = \frac{\log z}{\log_3 x}\exp\left(\sum_{k=1}^N \frac{(k-1)!C_k}{(\log_3 x)^k} \right)\left(1+O_N\left((\log_3 x)^{-(N+1)}\right)\right),\]
    where the numbers $C_k$ are the coefficients of the series expansion of $\Gamma(1+w)$ about $w=0$. 
    
    Therefore,
    \begin{align*}
    B(x) =& \frac{e^{-\gamma}x}{\log_3 x}\exp\left(\sum_{k=1}^N \frac{(k-1)!C_k}{(\log_3 x)^k} \right)\left(1+O_N\left((\log_3 x)^{-(N+1)}\right)\right) \\
    &\cdot \sum_{y<p\leq z} \frac{\exp(-\log_2 x / p)}{p^2}+O\left(\frac{x}{\log_2 x \exp(K\sqrt{\log_3 x})}\right)
    \end{align*}
    
    Our final task is to compute the sum over $p$. Write 
    \[\sum_{y<p\leq z} \frac{\exp(-\log_2 x / p)}{p^2} = \int_y^z \frac{\exp(-\log_2 x/t)}{t^2\log t} \, dt + O\left( \frac{x}{\exp(K\sqrt{\log_3 x})} \right)\]
    by again applying the prime number theorem with de la Vall\'ee Poussin error term. Make the change of variable $u=\log_2 x / t$, so that the integral becomes
    \begin{align*}
        \int_y^z \frac{\exp(-\log_2 x/t)}{t^2\log t} \, dt&=\frac{1}{\log_2 x} \int_{\log_2 x/z}^{\log_2 x / y} \frac{e^{-u}}{\log_3 x-\log u}\, du \\
        &= \frac{1}{\log_2 x \log_3 x} \int_{\log_2 x/z}^{\log_2 x / y} e^{-u} \left( 1-\frac{\log u}{\log_3 x} \right)^{-1} \, du.
    \end{align*}
    
    Now for any fixed $M$, we have
    \[\left( 1-\frac{\log u}{\log_3 x} \right)^{-1} = \sum_{k=0}^M \frac{\log^k u}{\log_3^k x} + O\left( \left(\frac{\log u}{\log_3 x}\right)^{M+1}  \right). \]
    Since $-\sqrt{\log_3 x} \leq \log u \leq \log_4 x$ in the domain of integration, the error term here is $O((\log_3 x)^{-(M+1)/2})$. We will also need to multiply this by $\int_{\log_2 x/z}^{\log_2 x/y}e^{-u}\,du \ll \log_3 x$.
    
    Recall the integral representation of the $\Gamma$ function,
    \[\Gamma(1+w)  =\int_0^{\infty} t^w e^{-t}\, dt, \]
    which converges absolutely for $\Re(w) > -1$. Using Leibniz' rule, the $k$th derivative of $\Gamma(1+w)$ equals
    \[\Gamma^{(k)}(1+w)= \int_0^{\infty} t^w e^{-t} \log^k t \, dt.\]
    We can now write
    \[\int_{\log_2 x/z}^{\log_2 x/y} e^{-u} \log^k u\, du = \Gamma^{(k)}(1) + O\left(\frac{x}{ \exp(K\sqrt{\log_3 x})} \right). \]
    
    Therefore, 
    \[\sum_{y<p\leq z}\frac{\exp(-\log_2 x/p)}{p^2} =\frac{1}{\log_2 x \log_3 x} \sum_{0\leq k \leq M} \frac{k!C_k}{(\log_3 x)^k} + O_M\left( \frac{1}{\log_2 x (\log_3 x)^{(M+1)/2}} \right).\]
    If we take $M=2N+3$, and absorb the extra terms appearing in the sum over $k$ into the error term, the error term is $O_N(1/(\log_2 x(\log_3 x)^{N+2})$.
    
   Assembling our results,
    
    \[B(x)=\frac{e^{-\gamma}x}{\log_2 x(\log_3 x)^2} \exp\left( \sum_{k=1}^N \frac{(k-1)!C_k}{(\log_3 x)^k} \right)\sum_{k=0}^N \frac{k! C_k}{(\log_3 x)^k}  + O\left( \frac{x}{\log_2 x (\log_3 x)^{N+3}} \right),\]
    which proves Theorem 1.
    
\section{Strictly nilpotent numbers: proof of Theorem 2}

The proof of Theorem 2 follows the same lines as the proof of Theorem 1.

Again let $y=\log_2 x/\log_3 x$, $z=\log_2 x\exp(\sqrt{\log_3 x})$, and $p$ be a prime in the interval $(y,z]$. Analogous to the proof of Theorem 1, define $S\p(x;p,0)$ to be the size of the set of $n\leq x$ such that $n=mp^3$, where $m$ is a product of distinct prime larger than $y$.

Define $S\p(x;p,k)$ to be the subset of $S\p(x;p,0)$ such that there are exactly $k$ primes dividing $m$ in the interval $(y,z]$, and at least one of these primes $q$ has the property that there is a prime $z<r\leq \log^{1/\log_2 x}$ dividing $m$ that is congruent to 1 modulo $q$ or 1 modulo $p$.

Define

\[D(x) = \sum_{y<p\leq z} \left( S\p(x;p,0)-\sum_{1\leq k\leq \log_3 x} S\p(x;p,k)\right).\]
We will show that $N(x)-A(x)=D(x)+O(x/z^2)$.

Suppose that $n$ is counted by $N(x)-A(x)$, but not $D(x)$. Then one of the following holds:
\begin{enumerate}
    \item There is a prime $q>z$ such that $q^3\mid n$.
    \item There are two distinct primes $q,r>y$ such that $q^3r^2$ divides $n$. 
    \item There is a prime $q$ in the interval $[2,y]$ dividing $n$.
    \item There is a prime $q>y$ such that $q^4\mid n$.
\end{enumerate}

The only condition that is not essentially the same as in the proof of Theorem 1 is (4), but the number of $n\leq x$ for which this holds is 
\[\sum_{q>y}\frac{x}{q^4} \ll \frac{x}{y^3\log y}\]
by Lemma \ref{primerec4}, which is $O(x/z^2)$.

Suppose that $n$ is counted by $D(x)$, but not $N(x)-A(x)$. Then one of the following holds:

\begin{enumerate}
    \item There is a prime $q>z$ such that $q$ divides $n$ and $\phi(n)$.
    \item $m$ has more than $\log_3 x$ prime factors in $(y,z]$.
    \item There is a prime $q$ in $(y,z]$ dividing $n$ such that there is a prime $r>x^{1/\log_2 x}$ dividing $m$ with $r$ congruent to 1 modulo $q$ or 1 modulo $p$.
    \item There is a prime $q$ in the interval $(y,z]$ dividing $m$ such that $p$ is congruent to 1 modulo $q$.
    \item There are two distinct primes $r$ and $s$ in $(y,z]$ dividing $n$ such that $r$ is congruent to 1 modulo $s$.
\end{enumerate}

Using arguments that are essentially the same as in the proof of Theorem 1, the number of $n$ such that one of these conditions hold is $O(x/z^2)$.

The estimation of $D(x)$ is identical to the estimation of $B(x)$, with the exception that we must compute the sum
\[\sum_{y<p\leq z} \frac{\exp(-\log_2 x/p)}{p^3} = \int_y^z \frac{\exp(-\log_2 x/t)}{t^3\log t}\, dt + O\left( \frac{x}{\exp(K\sqrt{\log_3 x})} \right),\]
again using the prime number theorem with de la Vall\'ee Poussin error term. Making the change of variables $u=\log_2 x/t$, we have for any fixed $M$

\[ \int_y^z \frac{\exp(-\log_2 x/t)}{t^3\log t}\, dt = \sum_{k=0}^M \frac{1}{(\log_2 x)^2 (\log_3 x)^{k+1}} \int_{\log_2 x/z}^{\log_2 x/y} u e^{-u} (\log u)^k \, du, \] 
where the error is $O((\log_2 x)^{-2}(\log_3 x)^{-(M+3)/2})$ (using that $\int_{\log_2 x/z}^{\log_2 x/y} ue^{-u}\, du \ll (\log_3 x)^2$).

Now \[\int_{\log_2 x/z}^{\log_2 x/y}ue^{-u}\log^k u\, du = \Gamma^{(k)}(2) + O\left( \frac{1}{\exp(K\sqrt{\log_3 x})} \right),\]
so that
\[\sum_{y<p\leq z}\frac{\exp(-\log_2 x/p)}{p^3} = \frac{1}{(\log_2 x)^2\log_3x}\sum_{k=0}^M\frac{\Gamma^{(k)}(2)}{(\log_3 x)^k} + O_M\left( \frac{1}{(\log_2 x)^2(\log_3 x)^{(M+3)/2}} \right). \]

Now setting $M=2N+1$ and absorbing the extra terms appearing in the sum over $k$ into the error term, we have
\[D(x)=\frac{e^{-\gamma x}}{(\log_2 x)^2(\log_3 x)^2}\exp\left(\sum_{k=1}^N \frac{(k-1)!C_k}{(\log_3 x)^k} \right)\sum_{k=0}^N \frac{\Gamma^{(k)}(2)}{(\log_3 x)^k} + O\left( \frac{x}{(\log_2 x)^2(\log_3 x)^{N+3}} \right).\]

    \section*{Acknowledgements}
    
    The author was partially supported by the Research and Training Group grant DMS-1344994 funded by the National Science Foundation. He thanks Paul Pollack for helpful comments.

\end{document}